\documentclass[a4paper,11pt]{amsart}
\usepackage{graphicx}
\usepackage{mathptmx}
\usepackage{amsmath}
\usepackage{amssymb}
\usepackage{enumitem}
\usepackage{xcolor}

\newmuskip\pFqmuskip

\newcommand*\pFq[6][8]{%
  \begingroup 
  \pFqmuskip=#1mu\relax
  \mathcode`=\string"8000
  \begingroup\lccode`\~=`\,
  \lowercase{\endgroup\let~}\pFqcomma
  F^{#2}_{#3}{\left(\genfrac..{0pt}{}{#4}{#5}\bigg|#6\right)}%
  \endgroup
}
\newcommand{\pFqcomma}{\mskip\pFqmuskip}

\newtheorem{theorem}{Theorem}
\newtheorem{lemma}[theorem]{Lemma}
\newtheorem{corollary}[theorem]{Corollary}
\newtheorem{proposition}[theorem]{Proposition}

\begin{document}

\title[Some identities involving derangement polynomials and numbers]{Some identities involving derangement polynomials and numbers and moments of gamma random variables}

\author{Taekyun  Kim}
\address{Department of Mathematics, Kwangwoon University, Seoul 139-701, Republic of Korea}
\email{tkkim@kw.ac.kr}

\author{Dae San Kim}
\address{Department of Mathematics, Sogang University, Seoul 121-742, Republic of Korea}
\email{dskim@sogang.ac.kr}

\author{Lee-Chae Jang}
\address{Graduate School of Eduction, Konkuk University, Seoul 143-701, Republic of Korea}
\email{lcjang@konkuk.ac.kr}

\author{Hyunseok Lee}
\address{Department of Mathematics, Kwangwoon University, Seoul 139-701, Republic of Korea}
\email{luciasconstant@kw.ac.kr}

\subjclass[2010]{11B73; 11B83; 65C50}
\keywords{derangement polynomials; cosine-derangement polynomials; sine-derangement polynomials; gamma random variable}

\begin{abstract}
The problem of counting derangements was initiated by Pierre R\'emonde de Motmort in 1708. A derangement is a permutation that has no fixed points and the derangement number $D_{n}$ is the number of fixed point free permutations on an $n$ element set. Furthermore, the derangement polynomials are natural extensions of the derangement numbers. In this paper, we study the derangement polynomials and numbers, their connections with cosine-derangement polynomials and sine-derangement polynomials and their applications to moments of some variants of gamma random variables.
\end{abstract}

\maketitle

\section{Introduction and preliminaries}

The problem of counting derangements was initiated by Pierre R\'emonde de Motmort in 1708 (see [1,2]).
A derangement is a permutation of the elements of a set, such that no element appears in its original position. In other words, a derangement is a permutation that has no fixed points. The derangement number $D_{n}$ is the number of fixed point free permutations on an $n\,(n \ge 1)$ element set. \par
The aim of this paper is to study derangement polynomials and numbers, their connections with cosine-derangement polynomials and sine-derangement polynomials and their applications to momnets of some variants of gamma random variables. Here the derangement polynomials $D_{n}(x)$ are natural extensions of the derangement numbers. \par
The outline of our main results is as follows. We show a recurrence relation for derangement polynomials. Then we derive identities involving derangement polynomials, Bell polynomials and Stirling numbers of both kinds. In addition, we also have an identity relating Bell polynomials, derangement polynomials and Euler numbers. Next, we introduce the two variable polynomials, namely cosine-derangement polynomials $D_n^{(c)}(x,y)$ and sine-derangement polynomials $D_n^{(s)}(x,y)$, in a natural manner by means of derangement polynomials. We obtain, among other things, their explicit expressions and recurrence relations. Lastly, in the final section we show that, if $X$ is the gamma random variable with parameters 1, 1, then $D_n(p), D_n^{(c)}(p,q), D_n^{(s)}(p,q)$ are given by the `moments' of some variants of $X$. \par
In the rest of this section, we recall the derangement numbers, especially their explicit expressions, generating function and recurrence relations. Also, we give the derangement polynomials and give their explicit expressions. Then we recall the gamma random variable with parameters $\alpha, \lambda$ along with their moments and the Bell polynomials. Finally, we give the definitions of the Stirling numbers of the first and second kinds. \par
As before, let $D_n $ denote the derangement number for $n \ge 1$, and let $D_0=1$. Then the first few derangement numbers $D_{n}\, (n \ge 0)$ are 1, 0, 1, 2, 9, 44, 265, 1854, 14833, 133496, 1334961, $\dots$. For $n\ge 0$, the derangement numbers are given by
\begin{align}
D_{n}\ &=\ n!-\binom{n}{1}(n-1)!+\binom{n}{2}(n-2)!-\binom{n}{3}(n-3)!+\cdots+(-1)^{n}\binom{n}{n}0!\label{1}\\
&\ =\ \sum_{k=0}^{n}\binom{n}{k}(n-k)!(-1)^{k}=n!\sum_{k=0}^{n}\frac{(-1)^{k}}{k!},\nonumber\quad(\mathrm{see}\ [8,10,11]).
\end{align}
From \eqref{1}, we note that
\begin{equation}
\sum_{n=0}^{\infty}D_{n}\frac{t^{n}}{n!}=\frac{1}{1-t}e^{-t},\quad(\mathrm{see}\ [1,2,3,8,10,14]).\label{2}
\end{equation}
By \eqref{2}, we get
\begin{equation}
e^{-t}=(1-t)\sum_{n=0}^{\infty}D_{n}\frac{t^{n}}{n!}=1+\sum_{n=1}^{\infty}\big(D_{n}-nD_{n-1}\big)\frac{t^{n}}{n!}.\label{3}	
\end{equation}
From \eqref{3}, we can easily derive the following recurrence relation :
\begin{equation}
(-1)^{n}=D_{n}-nD_{n-1},\quad(n\ge 1),\quad (\mathrm{see}\ [4,6,11,12,16]). \label{4}	
\end{equation}
Now, we consider the derangement polynomials which are given by
\begin{equation}
\frac{e^{-t}}{1-t}e^{xt}=\sum_{n=0}^{\infty}D_{n}(x)\frac{t^{n}}{n!},\quad(\mathrm{see}\ [12]).\label{5}
\end{equation}
From \eqref{5}, we have
\begin{equation}
\sum_{n=0}^{\infty}D_{n}(x)\frac{t^{n}}{n!}=\frac{1}{1-t}e^{-t}e^{xt}=\sum_{n=0}^{\infty}\bigg(\sum_{l=0}^{n}\binom{n}{l}D_{l}x^{n-l}\bigg)\frac{t^{n}}{n!}.\label{6}
\end{equation}
By comparing the coefficients on both sides of \eqref{6}, we get
\begin{equation}
D_{n}(x)=\sum_{l=0}^{n}\binom{n}{l}D_{l}x^{n-l},\quad (n\ge 0),\quad (\mathrm{see}\ [12]). \label{7}
\end{equation}
On the other hand,
\begin{align}
\frac{e^{-t}}{1-t}e^{xt}\ &=\ \frac{1}{1-t}e^{(x-1)t}\ =\ \sum_{l=0}^{\infty}t^{l}\sum_{m=0}^{\infty}(x-1)^{m}\frac{t^{m}}{m!}\label{8}\\
&=\ \sum_{n=0}^{\infty}\bigg(n!	\sum_{m=0}^{n}\frac{(x-1)^{m}}{m!}\bigg)\frac{t^{n}}{n!}.\nonumber
\end{align}
From \eqref{6},\eqref{7} and \eqref{8}, we have
\begin{equation}
D_{n}(x)=n!\sum_{m=0}^{n}\frac{(x-1)^{m}}{m!}=\sum_{l=0}^{n}\binom{n}{l}D_{l}x^{n-l},\quad(n\ge 0).\label{9}
\end{equation}
A continuous random variable $X$ whose density function is given by
\begin{equation}
f(x)=\left\{\begin{array}{ccc}
	\lambda e^{-\lambda x}(\lambda x)^{\alpha-1}/ \Gamma(\alpha), & \textrm{if $x\ge 0$}, \\
	0 & \textrm{if. $x<0$}
\end{array}\right.,\quad(\mathrm{see}\ [13,15,17]),\label{10}
\end{equation}
for some $\lambda>0$ and $\alpha>0$ is said to be the gamma random variable with parameter $\alpha,\lambda$ which is denoted by $X\sim\Gamma(\alpha,\lambda)$. \par
For $X\sim\Gamma(\alpha,\lambda)$, the $n$-th moment of $X$ is given by
\begin{align}
E[X^{n}]\ &=\ \frac{\lambda}{\Gamma(\alpha)}\int_{0}^{\infty}x^{n}e^{-\lambda x}(\lambda x)^{\alpha-1}dx. \label{11}\\
&=\ \frac{1}{\lambda^{n}\Gamma(\alpha)}	\int_{0}^{\infty}t^{n+\alpha-1}e^{-t}dt=\frac{\Gamma(\alpha+n)}{\lambda^{n}\Gamma(\alpha)}=\frac{(\alpha+n)\cdots(\alpha+1)\alpha}{\lambda^{n}}.\nonumber
\end{align}
It is well known that the Bell polynomials are defined by
\begin{equation}
e^{x(e^{t}-1)}=\sum_{n=0}^{\infty}\mathrm{Bel}_{n}(x)\frac{t^{n}}{n!},\quad(\mathrm{see}\ [9]).\label{12}
\end{equation}
When $x=1$, $\mathrm{Bel}_{n}=\mathrm{Bel}_{n}(1)$, $(n\ge 0)$ are called the Bell numbers. \par
The Stirling numbers of the first kind are defined as
\begin{equation}	
(x)_{n}=\sum_{l=0}^{n}S_{1}(n,l)x^{l},\quad(n\ge 0),\quad(\mathrm{see}\ [7,18]),\label{13}
\end{equation}
where $(x)_{0}=1,\ (x)_{n}=x(x-1)\cdots(x-n+1)$, $(n\ge 1)$. \par
As an inversion formula of \eqref{13}, the Stirling numbers of the second kind are defined by
\begin{equation}
x^{n}=\sum_{l=0}^{n}S_{2}(n,l)(x)_{l}\quad(n\ge 0),\quad(\mathrm{see}\ [5,7,18]).\label{14}
\end{equation}

\section{Derangement polynomials and numbers}
From \eqref{5}, we have
\begin{equation}
e^{(x-1)t}=\bigg(\sum_{n=0}^{\infty}D_{n}(x)\frac{t^{n}}{n!}\bigg)\big(1-t\big)=1+\sum_{n=1}^{\infty}\big(D_{n}(x)-nD_{n-1}(x)\big)\frac{t^{n}}{n!}.\label{15}
\end{equation}
On the other hand,
\begin{equation}
	e^{(x-1)t}=\sum_{n=0}^{\infty}\frac{(x-1)^{n}}{n!}t^{n}=1+\sum_{n=1}^{\infty}\frac{(x-1)^{n}}{n!}t^{n}. \label{16}
\end{equation}
Therefore, by \eqref{15} and \eqref{16}, we obtain the following lemma.
\begin{lemma}
	For $n\ge 1$, we have
	\begin{displaymath}
		D_{n}(x)-nD_{n-1}(x)=(x-1)^{n}.
	\end{displaymath}
\end{lemma}
Replacing $t$ by $1-e^{t}$ in \eqref{5}, we get
\begin{align}
e^{(1-x)(e^{t}-1)}\ &=\ e^{t}\sum_{l=0}^{\infty}D_{l}(x)\frac{1}{l!}(1-e^{t})^{l} \label{17}\\
&=\ \sum_{m=0}^{\infty}\frac{t^{m}}{m!}\sum_{l=0}^{\infty}(-1)^{l}D_{l}(x)\sum_{j=l}^{\infty}S_{2}(j,l)\frac{t^{j}}{j!}\nonumber\\
&=\ \sum_{m=0}^{\infty}\frac{t^{m}}{m!}\sum_{j=0}^{\infty}\bigg(\sum_{l=0}^{j}(-1)^{l}D_{l}(x)S_{2}(j,l)\frac{t^{j}}{j!}\nonumber\\
&=\ \sum_{n=0}^{\infty}\bigg(\sum_{j=0}^{n}\sum_{l=0}^{j}\binom{n}{j}(-1)^{l}D_{l}(x)S_{2}(j,l)\bigg)\frac{t^{n}}{n!}.\nonumber
\end{align}
From \eqref{17}, we have
\begin{equation}
\mathrm{Bel}_{n}(1-x)=\sum_{j=0}^{n}\sum_{l=0}^{j}\binom{n}{j}(-1)^{l}D_{l}(x)S_{2}(j,l),\quad(n\ge 0).\label{18}
\end{equation}
It is easy to show that
\begin{align}
\frac{1}{e^{t}}e^{(1-x)(e^{t}-1)}\ &=\ \sum_{l=0}^{\infty}\frac{(-1)^{l}}{l!}t^{l}\sum_{m=0}^{\infty}\mathrm{Bel}_{m}(1-x)\frac{t^{m}}{m!}\label{19} \\
&=\ \sum_{n=0}^{\infty}\bigg(\sum_{m=0}^{n}\binom{n}{m}\mathrm{Bel}_{m}(1-x)(-1)^{n-m}\bigg)\frac{t^{n}}{n!}.\nonumber	
\end{align}
Replacing $t$ by $\log(1-t)$ in \eqref{19}, we get
\begin{align}
\frac{1}{1-t}e^{-t}e^{xt}\ &=\ \sum_{l=0}^{\infty}\sum_{m=0}^{l}\binom{l}{m}\mathrm{Bel}_{m}(1-x)(-1)^{l-m}\frac{1}{l!}\big(\log(1-t)\big)^{l}\label{20} \\
&=\ \sum_{l=0}^{\infty}\sum_{m=0}^{l}\binom{l}{m}\mathrm{Bel}_{m}(1-x)(-1)^{l-m}\sum_{n=l}^{\infty}(-1)^{n}S_{1}(n,l)\frac{t^{n}}{n!}. \nonumber\\
&=\ \sum_{n=0}^{\infty}\bigg(\sum_{l=0}^{n}\sum_{m=0}^{l}\binom{l}{m}\mathrm{Bel}_{m}(1-x)(-1)^{n-l-m}S_{1}(n,l)\bigg)\frac{t^{n}}{n!}.\nonumber
\end{align}
From \eqref{5} and \eqref{20}, we have
\begin{equation}
D_{n}(x)=\sum_{l=0}^{n}\sum_{m=0}^{l}\binom{l}{m}\mathrm{Bel}_{m}(1-x)(-1)^{n-m-l}S_{1}(n,l),\quad(n\ge 0). \label{21}	
\end{equation}
Therefore, by \eqref{18} and \eqref{21}, we obtain the following theorem.
\begin{theorem}
	For $n\ge 0$, we have
	\begin{displaymath}
		\mathrm{Bel}_{n}(1-x)=\sum_{j=0}^{n}\sum_{l=0}^{j}\binom{n}{j}(-1)^{l}D_{l}(x)S_{2}(j,l),
	\end{displaymath}
	and
	\begin{displaymath}
		D_{n}(x)=\sum_{l=0}^{n}\sum_{m=0}^{l}\binom{l}{m}\mathrm{Bel}_{m}(1-x)(-1)^{n-m-l}S_{1}(n,l).
	\end{displaymath}
\end{theorem}
\begin{corollary}
	For $n\ge 0$, we have
		\begin{displaymath}
		\mathrm{Bel}_{n}=\sum_{j=0}^{n}\sum_{l=0}^{j}\binom{n}{j}(-1)^{l}D_{l}S_{2}(j,l),
	\end{displaymath}
	and
	\begin{displaymath}
		D_{n}=\sum_{l=0}^{n}\sum_{m=0}^{l}\binom{l}{m}\mathrm{Bel}_{m}(-1)^{n-m-l}S_{1}(n,l).
	\end{displaymath}
\end{corollary}
Replacing $t$ by $-e^{t}$ in \eqref{5}, we get
\begin{align}
\frac{1}{e^{t}+1}e^{(1-x)e^{t}}\ &=\ \sum_{m=0}^{\infty}D_{m}(x)\frac{(-1)^{m}}{m!}e^{mt} \nonumber\\
&=\ \sum_{m=0}^{\infty}\frac{D_{m}(x)(-1)^{m}}{m!}\sum_{n=0}^{\infty}m^{n}\frac{t^{n}}{n!}\label{22}\\
&=\ \sum_{n=0}^{\infty}\bigg(\sum_{m=0}^{\infty}\frac{(-1)^{m}D_{m}(x)}{m!}m^{n}\bigg)\frac{t^{n}}{n!}.\nonumber
\end{align}
On the other hand, we have
\begin{align}
\frac{1}{e^{t}+1}e^{(1-x)e^{t}}\ &=\ \frac{e^{1-x}}{2}\frac{2}{e^{t}+1}e^{(1-x)(e^{t}-1)} \label{23}\\
&=\ \frac{e^{1-x}}{2}\sum_{l=0}^{\infty}E_{l}\frac{t^{l}}{l!}\sum_{m=0}^{\infty}\mathrm{Bel}_{m}(1-x)\frac{t^{m}}{m!}\nonumber\\
&=\ \frac{e^{1-x}}{2}\sum_{n=0}^{\infty}\bigg(\sum_{m=0}^{n}\mathrm{Bel}_{m}(1-x)E_{n-m}\binom{n}{m}\bigg)\frac{t^{n}}{n!},\nonumber
	\end{align}
where $E_{n}$ are the ordinary Euler numbers. \par
Therefore, by \eqref{22} and \eqref{23}, we obtain the following theorem,
\begin{theorem}
	For $n\ge 0$, we have
	\begin{displaymath}
		\sum_{m=0}^{n}\mathrm{Bel}_{m}(1-x)E_{n-m}\binom{n}{m}=2e^{x-1}\sum_{m=0}^{\infty}(-1)^{m}\frac{D_{m}(x)}{m!}m^{n}.
	\end{displaymath}
\end{theorem}
Now, we observe that
\begin{align}
\bigg(\frac{1}{1-t}\bigg)^{r}\ &=\ \bigg(\frac{1}{1-t}\bigg)^{r}e^{-rt}e^{rt}=\bigg(\frac{1}{1-t}e^{-t}\bigg)^{r-1}\frac{e^{-t}}{1-t}e^{rt}	\nonumber \\
&=\ \sum_{k=0}^{\infty}\sum_{l_{1}+\cdots+l_{r-1}=k}\binom{k}{l_{1},\dots,l_{r-1}}D_{l_{1}}D_{l_{2}}\cdots D_{l_{r}}\frac{t^{k}}{k!}\sum_{m=0}^{\infty}D_{m}(r)\frac{t^{m}}{m!}\label{24} \\
&=\ \sum_{n=0}^{\infty}\bigg(\sum_{k=0}^{n}\sum_{l_{1}+\cdots+l_{r-1}=k}\binom{k}{l_{1},\dots,l_{r-1}}\binom{n}{k}D_{l_{1}}D_{l_{2}}\cdots D_{l_{k-1}}D_{n-k}(r)\bigg)\frac{t^{n}}{n!},\nonumber
\end{align}
where $r$ is a positive integer. \par
On the other hand,
\begin{equation}
\bigg(\frac{1}{1-t}\bigg)^{r}=\sum_{n=0}^{\infty}\binom{-r}{n}(-1)^{n}t^{n}=\sum_{n=0}^{\infty}n!\binom{r+n-1}{n}\frac{t^{n}}{n!}. \label{25}
\end{equation}
Therefore, by \eqref{24} and \eqref{25}, we obtain the following Proposition.
\begin{proposition}
For $r\in\mathbb{N}$, we have
\begin{displaymath}
\binom{r+n-1}{n}=\frac{1}{n!}\sum_{k=0}^{n}\sum_{l_{1}+\cdots+l_{r-1}=k}\binom{k}{l_{1},\dots,l_{r-1}}\binom{n}{k}D_{l_{1}}\cdots D_{l_{r-1}}D_{n-k}(r).
\end{displaymath}
\end{proposition}
It is well known that
\begin{equation}
e^{ix}=\cos x+i\sin x,\quad i=\sqrt{-1},\quad(\mathrm{see}\ [5,7,19]).\label{26}
\end{equation}
From \eqref{5}, we note that
\begin{equation}
\frac{e^{-t}}{1-t}e^{(x+iy)t}=\sum_{n=0}^{\infty}D_{n}(x+iy)\frac{t^{n}}{n!},\quad(x,y\in\mathbb{R}),\label{27}
\end{equation}
and
\begin{equation}
\frac{e^{-t}}{1-t}e^{(x-iy)t}=\sum_{n=0}^{\infty}D_{n}(x-iy)\frac{t^{n}}{n!}.\label{28}	
\end{equation}
By \eqref{9}, \eqref{27} and \eqref{28}, we get
\begin{equation}
D_{n}(x+iy)=n!\sum_{m=0}^{n}\frac{(x-1+iy)^{m}}{m!},\label{29}
\end{equation}
and
\begin{equation}
D_{n}(x-iy)=n!\sum_{m=0}^{n}\frac{(x-1-iy)^{m}}{m!},\quad(n\ge 0).\label{30}
\end{equation}
From \eqref{29} and \eqref{30}, we can derive the following equations:
\begin{equation}
\frac{e^{-t}}{1-t}e^{xt}\cos(yt)=\sum_{n=0}^{\infty}\bigg(\frac{D_{n}(x+iy)+D_{n}(x-iy)}{2}\bigg)\frac{t^{n}}{n!},\label{31}	
\end{equation}
and
\begin{equation}
\frac{e^{-t}}{1-t}e^{xt}\sin(yt)=\sum_{n=0}^{\infty}\bigg(\frac{D_{n}(x+iy)-D_{n}(x-iy)}{2i}\bigg)\frac{t^{n}}{n!}.\label{32}
\end{equation}
We define cosine-derangement polynomials and sine-derangement polynomials respectively by
\begin{equation}
\frac{e^{-t}}{1-t}e^{xt}\cos yt	=\sum_{n=0}^{\infty}D_{n}^{(c)}(x,y)\frac{t^{n}}{n!},\label{33}
\end{equation}
and
\begin{equation}
\frac{e^{-t}}{1-t}e^{xt}\sin yt	=\sum_{n=0}^{\infty}D_{n}^{(s)}(x,y)\frac{t^{n}}{n!},\label{34}
\end{equation}
Thus, we have
\begin{equation}
D_{n}^{(c)}(x,y)=\frac{D_{n}(x+iy)+D_{n}(x-iy)}{2},\label{35}	
\end{equation}
and
\begin{equation}
	D_{n}^{(s)}(x,y)=\frac{D_{n}(x+iy)-D_{n}(x-iy)}{2i},\quad(n\ge 0).\label{36}
\end{equation}
Therefore, we obtain the following theorem.
\begin{theorem}
	For $n\ge 0$, we have
	\begin{displaymath}
		D_{n}^{(c)}(x,y)=\frac{n!}{2}\sum_{m=0}^{n}\frac{1}{m!}\big((x-1+iy)^{m}+(x-1-iy)^{m}\big),
	\end{displaymath}
	and
	\begin{displaymath}
		D_{n}^{(s)}(x,y)=\frac{n!}{2i}\sum_{m=0}^{n}\frac{1}{m!}\big((x-1+iy)^{m}-(x-1-iy)^{m}\big).
	\end{displaymath}	
\end{theorem}
From \eqref{33}, we note that
\begin{align}
\sum_{n=0}^{\infty}D_{n}^{(c)}(x,y)\frac{t^{n}}{n!}\ &=\ \frac{e^{-t}}{1-t}e^{xt}\cos(yt)\label{37}\\
&=\ \sum_{l=0}^{\infty}\frac{D_{l}}{l!}t^{l}\sum_{k=0}^{\infty}\sum_{m=0}^{[\frac{k}{2}]}\binom{k}{2m}(-1)^{m}y^{2m}x^{k-2m}\frac{t^{k}}{k!}	\nonumber \\
&=\ \sum_{n=0}^{\infty}\bigg(\sum_{k=0}^{n}\binom{n}{k}D_{n-k}\sum_{m=0}^{[\frac{k}{2}]}\binom{k}{2m}(-1)^{m}y^{2m}x^{k-2m}\bigg)\frac{t^{n}}{n!}.\nonumber
\end{align}
Therefore, by comparing the coefficients on both sides of \eqref{37}, we obtain the following theorem.
\begin{theorem}
	For $n\ge 0$, we have
	\begin{displaymath}
		D_{n}^{(c)}(x,y)=\sum_{m=0}^{[\frac{n}{2}]}\sum_{k=2m}^{n}\binom{n}{k}\binom{k}{2m}D_{n-k}(-1)^{m}y^{2m}x^{k-2m}.
	\end{displaymath}
\end{theorem}
\begin{corollary}
	For $n\ge 0$, we have
	\begin{displaymath}
		\frac{n!}{2}\sum_{m=0}^{n}\frac{1}{m!}\big((x-1+iy)^{m}+(x-1-iy)^{m}\big)=\sum_{m=0}^{[\frac{n}{2}]}\sum_{k=2m}^{n}\binom{n}{k}\binom{k}{2m}D_{n-k}(-1)^{m}y^{2m}x^{k-2m}.
	\end{displaymath}
\end{corollary}
By \eqref{33}, we get
\begin{align}
e^{(x-1)t}\cos yt\ &=\ (1-t)\sum_{n=0}^{\infty}D_{n}^{(c)}(x,y)\frac{t^{n}}{n!} \label{38} \\
&=\ 1+\sum_{n=1}^{\infty}\big(D_{n}^{(c)}(x,y)-nD_{n-1}^{(c)}(x,y)\big)\frac{t^{n}}{n!}.\nonumber
\end{align}
Thus, we have
\begin{align}
\cos yt\ &=\ e^{(1-x)t}+e^{(1-x)t}\sum_{m=1}^{\infty}\big(D_{m}^{(c)}(x,y)-mD_{m-1}^{(c)}(x,y)\big)\frac{t^{m}}{m!} \label{39} \\
&=\ \sum_{n=0}^{\infty}(1-x)^{n}\frac{t^n}{n!}+\sum_{l=0}^{\infty}(1-x)^{l}\frac{t^l}{l!}\sum_{m=1}^{\infty}\big(D_{m}^{(c)}(x,y)-mD_{m-1}^{(c)}(x,y)\big)\frac{t^{m}}{m!}\nonumber\\
&=\ 1+\sum_{n=1}^{\infty}\bigg((1-x)^{n}+\sum_{m=1}^{n}\binom{n}{m}(1-x)^{n-m}\big(D_{m}^{(c)}(x,y)-mD_{m-1}^{(c)}(x,y)\big)\bigg)\frac{t^{n}}{n!}.\nonumber
\end{align}
On the other hand, we also have
\begin{equation}
\cos yt=\sum_{n=0}^{\infty}\frac{(-1)^{n}}{(2n)!}y^{2n}t^{2n}. \label{40}	
\end{equation}
Therefore, by \eqref{39} and \eqref{40}, we obtain the following theorem.
\begin{theorem}
	For $k\in\mathbb{N}$, we have
	\begin{displaymath}
		(1-x)^{n}+\sum_{m=1}^{n}\binom{n}{m}(1-x)^{n-m}\big(D_{m}^{(c)}(x,y)-mD_{m-1}^{(c)}(x,y)\big)=\left\{\begin{array}{ccc}
			(-1)^{k}y^{2k}, & \textrm{if $n=2k$,}\\
			0, & \textrm{if $n=2k-1$}.
		\end{array}\right.
	\end{displaymath}
\end{theorem}
By \eqref{33}, we get
\begin{align}
e^{(x-1)t}\cos yt\ &=\ \sum_{n=0}^{\infty}D_{n}^{(c)}(x,y)\frac{t^{n}}{n!}(1-t) \label{41}\\
&=\ \sum_{n=1}^{\infty}\big(D_{n}^{(c)}(x,y)-nD_{n-1}^{(c)}(x,y)\big)\frac{t^{n}}{n!}+1.\nonumber
\end{align}
On the other hand,
\begin{align}
e^{(x-1)t}\cos yt\ &=\ \sum_{l=0}^{\infty}(x-1)^{l}\frac{t^{l}}{l!}\sum_{m=0}^{\infty}y^{2m}(-1)^{m}\frac{t^{2m}}{(2m)!} \label{42} \\
&=\ 1+\sum_{n=1}^{\infty}\bigg(\sum_{m=0}^{[\frac{n}{2}]}\binom{n}{2m}(-1)^m(x-1)^{n-2m}y^{2m}\bigg)\frac{t^{n}}{n!}.\nonumber
\end{align}
Therefore, by \eqref{41} and \eqref{42}, we obtain the following theorem.
\begin{theorem}
	For $n\ge 1$, we have
	\begin{displaymath}
		D_{n}^{(c)}(x,y)-nD_{n-1}^{(c)}(x,y)=\sum_{m=0}^{[\frac{n}{2}]}\binom{n}{2m}(-1)^{m}(x-1)^{n-2m}y^{2m}.
	\end{displaymath}
\end{theorem}
It is not difficult to show that
\begin{equation}
\sum_{n=0}^{\infty}D_{n}^{(c)}(x+r,y)\frac{t^{n}}{n!}=\sum_{n=0}^{\infty}\bigg(\sum_{l=0}^{n}\binom{n}{l}D_{l}^{(c)}(x,y)r^{n-l}\bigg)\frac{t^{n}}{n!},\label{43}
\end{equation}
where $r$ is positive integer. \par
By comparing the coefficients on both sides of \eqref{39}, we get
\begin{equation}
D_{n}^{(c)}(x+r,y)=\sum_{l=0}^{n}\binom{n}{l}D_{l}^{(c)}(x,y)r^{n-l}. \label{44}	
\end{equation}
Now, we observe that
\begin{equation}
\sum_{n=1}^{\infty}\frac{\partial}{\partial x}D_{n}^{(c)}(x,y)\frac{t^{n}}{n!}=\frac{\partial}{\partial x}\bigg(\frac{e^{-t}}{1-t}e^{xt}\cos yt \bigg)\label{45}\qquad\qquad\qquad
\end{equation}
\begin{displaymath}
	=t\frac{e^{-t}}{1-t}e^{xt}\cos yt=t\sum_{n=0}^{\infty}D_{n}^{(c)}(x,y)\frac{t^{n}}{n!}=\sum_{n=1}^{\infty}nD_{n-1}^{(c)}(x,y)\frac{t^{n}}{n!}.
\end{displaymath}
Form \eqref{45}, we note that
\begin{displaymath}
	D_{0}^{(c)}(x,y)=1,\quad\frac{\partial}{\partial x}D_{n}^{(c)}(x,y)=nD_{n-1}^{(c)}(x,y),\quad(n\ge 1).
\end{displaymath}
Therefore, we obtain the following theorem.
\begin{theorem}
	For $n\ge 0$, we have
	\begin{displaymath}
		D_{0}^{(c)}(x,y)=0,\quad \frac{\partial}{\partial x}D_{n}^{(c)}(x,y)=nD_{n-1}^{(c)}(x,y),\quad(n\ge 1).
	\end{displaymath}
	In particular,
	\begin{displaymath}
		\frac{d}{dx}D_{n}(x)=\frac{\partial}{\partial x}D_{n}^{(c)}(x,0)=nD_{n-1}^{(c)}(x,0)=nD_{n-1}^{(c)}(x),\quad (n\ge 1).
	\end{displaymath}
\end{theorem}
\begin{corollary}
$D_{n}^{(c)}(x,y)$ as a polynomial in $x$, for each fixed $y$, and $D_{n}(x)$ are Appell sequences.
\end{corollary}
From \eqref{34}, we note that
\begin{align}
\sum_{n=0}^{\infty}D_{n}^{(s)}(x,y)\frac{t^{n}}{n!}\ &=\ \frac{1}{1-t}e^{-t}e^{xt}\sin yt \label{46} \\
&=\ \sum_{k=0}^{\infty}\frac{D_{k}}{k!}t^{k}\sum_{j=1}^{\infty}\sum_{m=0}^{[\frac{j-1}{2}]}\binom{j}{2m+1}x^{j-2m-1}y^{2m+1}\frac{t^{j}}{j!}\nonumber \\
&=\ \sum_{n=1}^{\infty}\bigg(\sum_{j=1}^{n}\sum_{m=0}^{[\frac{j-1}{2}]}\binom{j}{2m+1}\binom{n}{j}x^{j-2m-1}y^{2m+1}D_{n-j}\bigg)\frac{t^{n}}{n!}.\nonumber
\end{align}
Therefore, by \eqref{46}, we obtain the following theorem.
\begin{theorem}
	For $n\ge 0$, we have
	\begin{displaymath}
		D_{0}^{(s)}(x,y)=0,\quad D_{n}^{(s)}(x,y)=\sum_{j=1}^{n}\sum_{m=0}^{[\frac{j-1}{2}]}\binom{j}{2m+1}\binom{n}{j}x^{j-2m-1}y^{2m+1}D_{n-j}.
	\end{displaymath}
\end{theorem}
By \eqref{30},\eqref{32} and Theorem 13, we obtain the following theorem.
\begin{corollary}
	For $n\ge 1$, we have
	\begin{displaymath}
		\frac{D_{n}(x+iy)-D_{n}(x-iy)}{2i}=\sum_{j=1}^{n}\sum_{m=1}^{[\frac{j-1}{2}]}\binom{j}{2m+1}\binom{n}{j}x^{j-2m-1}y^{2m+1}D_{n-j}.
	\end{displaymath}
\end{corollary}
By \eqref{46}, we see that
\begin{align}
\sin yt\ &=\ e^{(1-x)t}\sum_{k=1}^{\infty}\big(D_{k}^{(s)}(x,y)-kD_{k-1}^{(s)}(x,y)\big)\frac{t^{k}}{k!}\label{47}\\
&=\ 	\sum_{m=0}^{\infty}(1-x)^{m}\frac{t^m}{m!}\sum_{k=1}^{\infty}\big(D_{k}^{(s)}(x,y)-kD_{k-1}^{(s)}(x,y)\big)\frac{t^{k}}{k!}\nonumber \\
&=\ \sum_{n=1}^{\infty}\bigg(\sum_{k=1}^{n}\binom{n}{k}\big(D_{k}^{(s)}(x,y)-kD_{k-1}^{(s)}(x,y)\big)(1-x)^{n-k}\bigg)\frac{t^{n}}{n!}. \nonumber
\end{align}
On the other hand, we also have
\begin{equation}
\sin yt=\sum_{n=1}^{\infty}\frac{(-1)^{n-1}}{(2n-1)!}y^{2n-1}t^{2n-1}. \label{48}	
\end{equation}
Therefore, by \eqref{47} and \eqref{48}, we obtain the following theorem.
\begin{theorem}
For $m\in\mathbb{N}$, we have
\begin{displaymath}
\sum_{k=1}^{n}\binom{n}{k}\big(D_{k}^{(s)}(x,y)-kD_{k-1}^{(s)}(x,y)\big)(1-x)^{n-k}=\left\{\begin{array}
{ccc}
(-1)^{m-1}y^{2m-1}, & \textrm{if $n=2m-1,$}\\
0, & \textrm{if $n=2m.$\quad\ \ }
\end{array}\right..
\end{displaymath}
\end{theorem}
It is easy to show that $\displaystyle\frac{\partial}{\partial x}D_{n}^{(s)}(x,y)=nD_{n-1}^{(s)}(x,y)\displaystyle$. However, $D_{n}^{(s)}(x,y)$ is not an Appell sequence, since  $D_{0}^{(s)}(x,y)=0$. \par
We observe that
\begin{align}
\sum_{n=0}^{\infty}D_{n}^{(s)}(x,y)\frac{t^{n}}{n!}\ &=\ \frac{e^{-t}}{1-t}e^{xt}\sin yt \label{49}\\
&=\ \sum_{l=0}^{\infty}D_{l}(x)\frac{t^{l}}{l!}\sum_{m=0}^{\infty}(-1)^{m}y^{2m+1}\frac{t^{2m+1}}{(2m+1)!}\nonumber \\
&=\ \sum_{n=1}^{\infty}\bigg(\sum_{m=0}^{[\frac{n-1}{2}]}\binom{n}{2m+1}(-1)^{m}y^{2m+1}D_{n-2m-1}(x)\bigg)\frac{t^{n}}{n!}. \nonumber
\end{align}
Comparing the coefficients on both sides of \eqref{49}, we have the following theorem.
\begin{theorem}
For $n\ge 1$, we have
\begin{displaymath}
D_{n}^{(s)}(x,y)= \sum_{m=0}^{[\frac{n-1}{2}]}\binom{n}{2m+1}(-1)^{m}y^{2m+1}D_{n-2m-1}(x).
\end{displaymath}
\end{theorem}
For $r\in\mathbb{N}$, we have
\begin{align*}
\sum_{n=0}^{\infty}D_{n}^{(s)}(x+r,y)\ &=\ \frac{e^{-t}}{1-t}e^{(x+r)t}\sin yt\
=\ \frac{e^{-t}}{1-t}e^{xt}\sin yte^{rt}\\
&=\ \sum_{l=0}^{\infty}D_{l}^{(s)}(x,y)\frac{t^{l}}{l!}\sum_{m=0}^{\infty}r^{m}\frac{t^{m}}{m!}\\
&=\ \sum_{n=0}^{\infty}\bigg(\sum_{l=0}^{n}\binom{n}{l}D_{l}^{(s)}(x,y)r^{n-l}\bigg)\frac{t^{n}}{n!}.
\end{align*}
Thus, we obtain
\begin{displaymath}
D_{n}^{(s)}(x+r,y)=\sum_{l=0}^{n}\binom{n}{l}D_{l}^{(s)}(x,y)r^{n-l},\quad(n\ge 0).
\end{displaymath}

\section{Further Remarks}
Let $X$ be a gamma random variable with parameters 1,1 which is denoted by $X\sim\Gamma(1,1)$. \\
Then we observe that
\begin{equation}
E\big[e^{(X-1+p)t}\big]=\int_{0}^{\infty}e^{(x-1+p)t}f(x)dx,\label{50}	
\end{equation}
where $f(x)$ is the density function of $X$, and $p \in \mathbb{R}$. \par
From \eqref{10} and \eqref{50}, we can derive the following equation \eqref{51}.
\begin{align}
E\big[e^{(X-1+p)t}\big]\ &=\ \int_{0}^{\infty}e^{(x-1+p)t}e^{-x}dx \label{51} \\
&=\ e^{-t+pt}\cdot\int_{0}^{\infty}e^{-x(1-t)}dx\nonumber\\
&=\ \frac{e^{-t}}{1-t}e^{pt}\ =\ \sum_{n=0}^{\infty}D_{n}(p)\frac{t^{n}}{n!}.\nonumber	
\end{align}
On the other hand, by Taylor expansion, we get
\begin{equation}
E\big[e^{(X-1+p)t}\big]=\sum_{n=0}^{\infty}E\big[(X-1+p)^{n}\big]\frac{t^{n}}{n!}.\label{52}
\end{equation}
Therefore, by \eqref{51} and \eqref{52}, we obtain the following theorem.
\begin{theorem}
For $n\ge 0$, $X\sim\Gamma(1,1)$, the moment of $X-1+p$ is given by
\begin{displaymath}
E\big[(X-1+p)^{n}\big]=D_{n}(p).
\end{displaymath}
\end{theorem}
When $p=0$, $D_{n}=D_{n}(0)=E[(X-1)^{n}]$, $(n\ge 0)$. \par
Thus, we note that
\begin{equation}
D_{n}=\sum_{l=0}^{n}\binom{n}{l}(-1)^{n-l}E\big[X^{l}\big].\label{53}	
\end{equation}
For $X\sim\Gamma(1,1)$, we note that the moment of $X$ is given by $E[X^{n}]=n!$, $(n\ge 0)$. \par
Therefore, by \eqref{53}, we obtain the following corollary.
\begin{corollary}
	For $n\ge 0$, $X\sim\Gamma(1,1)$, we have
	\begin{displaymath}
		D_{n}=\sum_{l=0}^{n}\binom{n}{l}(-1)^{n-l}l!
	\end{displaymath}
and
\begin{displaymath}
	D_{n}(p)=\sum_{l=0}^{n}\binom{n}{l}(p-1)^{n-l}l!.
\end{displaymath}
\end{corollary}
For $X\sim\Gamma(1,1)$, we have
\begin{equation}
E\big[e^{(X-1+p+iq)t}\big]=\frac{e^{-t}}{1-t}e^{(p+iq)t},\quad\mathrm{where}\ p,q\in\mathbb{R}.\label{54}
\end{equation}
From \eqref{54}, we note that
\begin{equation}
	E\big[e^{(X-1+p-iq)t}\big]=\frac{e^{-t}}{1-t}e^{(p-iq)t}.\label{55}
\end{equation}
By \eqref{54} and \eqref{55}, we get
\begin{equation}
	E\big[e^{(X-1+p+iq)t}\big]+ E\big[e^{(X-1+p-iq)t}\big]=\frac{2e^{-t}}{1-t}e^{pt}\cos qt=\sum_{n=0}^{\infty}2D_{n}^{(c)}(p,q)\frac{t^{n}}{n!}. \label{56}
\end{equation}
On the other hand, by Taylor expansion, we get
\begin{align}
& E\big[e^{(X-1+p+iq)t}\big]+ E\big[e^{(X-1+p-iq)t}\big]\label{57} \\
&\quad=\ \sum_{n=0}^{\infty}E\big[(X-1+p+iq)^{n}+(X-1+p-iq)^{n}\big]\frac{t^{n}}{n!}.\nonumber
\end{align}
Therefore, by \eqref{56} and \eqref{57}, we obtain the following theorem.
\begin{theorem}
	For $n\ge 0$, $X\sim\Gamma(1,1)$, we have
	\begin{displaymath}
		E\bigg[\frac{(X-1+p+iq)^{n}+(X-1+p-iq)^{n}}{2}\bigg]=D_{n}^{(c)}(p,q).
	\end{displaymath}
\end{theorem}
It is easy to show that
\begin{align*}
	& E\big[e^{(X-1+p+iq)t}\big]-E\big[e^{(X-1+p-iq)t}\big]\\ \\
&\quad=\ 2i\frac{e^{-t}}{1-t}e^{pt}\sin qt=(2i)\sum_{n=1}^{\infty}D_{n}^{(s)}(p,q)\frac{t^{n}}{n!},
\end{align*}
where $X\sim\Gamma(1,1)$. \par
Thus, we have
\begin{displaymath}
	E\bigg[\frac{(X-1+p+iq)^{n}-(X-1+p-iq)^{n}}{2i}\bigg]=D_{n}^{(s)}(p,q),\quad(n\ge 0),
\end{displaymath}
where $X\sim\Gamma(1,1)$. \par

\section{Conclusion}
The introduction of deragement numbers $D_n$ goes back to as early as 1708 when
Pierre R\'emond de Montmort considered some counting problem on derangements. In this paper, we dealt with derangement polynomials $D_n(x)$ which are natural extensions of the derangement numbers. We showed a recurrence relation for derangement polynomials. We derived identities involving derangement polynomials, Bell polynomials and Stirling numbers of both kinds. In addition, we also obtained an identity relating Bell polynomials, derangement polynomials and Euler numbers. Next, we introduced the cosine-derangement polynomials $D_n^{(c)}(x,y)$ and sine-derangement polynomials $D_n^{(s)}(x,y)$, by means of derangement polynomials. Then we derived, among other things, their explicit expressions and recurrence relations. Lastly, as an applications we showed that, if $X$ is the gamma random variable with parameters 1, 1, then $D_n(p), D_n^{(c)}(p,q), D_n^{(s)}(p,q)$ are given by the `moments' of some variants of $X$. \par
We have witnessed that the study of some special numbers and polynomials was done intensively by using several different means, which include generating functions, combinatorial methods, umbral calculus, $p$-adic analysis, probability theory, special functions and differential equations. Moreover, the same has been done for various degenerate versions of quite a few special numbers and polynomials in recent years with their interests not only in combinatorial and arithmetical properties but also in their applications to symmetric identities, differential equations and probability theories. It would have been nicer if we were able to find abundant applications in other disciplines. \par
It is one of our future projects to continue to investigate many ordinary and degenerate special numbers and polynomials by various means and find their applications in physics, science and engineering as well as in mathematics.

\vspace{0.5cm}

{\bf Author Contributions:} T.K. and D.S.K. conceived of the framework and structured the whole paper; D.S.K. and
T.K. wrote the paper; L.C.J. checked the errors of the paper; H.L. typed the paper; D.S.K. and T.K.
completed the revision of the article. All authors have read and agreed to the published version of the manuscript.

\vspace{0.3cm}

{\bf Funding:} Not applicable.

\vspace{0.3cm}

{\bf Conflicts of Interest:} The authors declare no conflicts of interest.

\end{document}